\documentclass[11pt]{amsart}

\usepackage{amsmath}
\usepackage{amsthm}
\usepackage{amsfonts}

\newtheorem{defi}{\bf Definition}[section]
\newtheorem{teor}{\bf Theorem}

\newtheorem{propo}{\bf Proposition}[section]

\newtheorem{teo}{\bf Theorem}[section]
\newtheorem{lema}{\bf Lemma}[section]

\newtheorem{cor}{\bf Corollary}[section]

\begin{document}
\title[Quasi Anosov Diffeomorphisms]{On Manifolds Supporting Quasi Anosov Diffeomorphisms}
\author{Jana Rodriguez Hertz}
\address{Jana Rodriguez Hertz}
\email{jana@fing.edu.uy}
\thanks{The first author was partially supported by a grant from PEDECIBA}

\author{Ra\'{u}l
Ures} \address{Ra\'{u}l Ures} \email{ures@fing.edu.uy}
\thanks{The second author was partially supported by CONICYT, Fondo Clemente
Estable}

\author{Jos\'{e} L. Vieitez}
\address{Jos\'{e} L. Vieitez}
\email{jlvb@fing.edu.uy} \address{CC 30, IMERL - Facultad  de
Ingenier\'{\i}a
\\Universidad de la Rep\'{u}blica\\Montevideo, Uruguay}


\begin{abstract}
Let $M$ be an $n$-dimensional manifold supporting a quasi Anosov
diffeomorphism. If $n=3$ then either $M={\mathbb T}^3$, in which
case the diffeomorphisms is Anosov, or else its fundamental group
contains a copy of ${\mathbb Z} ^6$. If $n=4$ then $\Pi_1(M)$
contains a copy of ${\mathbb Z} ^4$, provided that the
diffeomorphism is not Anosov.
\end{abstract}

\maketitle
\thispagestyle{empty}
\section{Introduction}

In this work we obtain some restrictions on a manifold $M$ in
order to support {\em quasi Anosov diffeomorphisms} (QAD), meaning
diffeomorphisms $f:M\rightarrow M$ such that
$\|(f^n)'(x)v\|\rightarrow+\infty$ for all non zero vectors $v\in
T_xM$, either for $n\rightarrow\infty$ or for
$n\rightarrow-\infty$. \par These maps were introduced in
\cite{m1} by Ma\~{n}\'{e}, who showed they satisfy Axiom A and the
no-cycle condition. Besides, they are the $C^1$ interior of the
class of
 {\em expansive diffeomorphisms} $g$ \cite{m2}, those satisfying
for some $\alpha >0$, that $\sup_{n\in {\mathbb Z
}}d(f^n(x),f^n(y))>\alpha$, if $x\not=y$.

We expect on one hand, that studying the relation between QAD's
(seen as simplified examples of expansive diffeomorphisms) and the
manifolds supporting them, will bring some ideas about the
restrictions that expansive homeomorphisms impose on their ambient
manifolds. Indeed, let us recall, for instance, that no expansive
homeomorphism may be found on the 2-sphere, all expansive
homeomorphisms are conjugated to Anosov diffeomorphisms on the
2-torus, and to pseudo-Anosov maps if they live in a surface with
genus greater than 1 \cite{h,l}.

In the case $M$ is a 3-manifold, an expansive $C^{1+\varepsilon}$
diffeomorphism is conjugated to an Anosov diffeomorphism if its
non-wandering set is $M$; in which case $M$ must be $\mathbb
{T}^3$ \cite{v}.

On the other hand, Ma\~{n}\'{e} \cite{m3} and Hiraide \cite{h2} have posed
the question of whether a QAD on $\mathbb {T}^n$ is necessarily
Anosov. This work provides a positive answer for $n=3$, though we
must point out that solving the problem for higher dimensions
seems to involve more sophisticated tools.

However, there are also examples of QAD's on 3-dimensional
manifolds having wandering points \cite{fr}. We present here
necessary conditions on a 3-manifold to support this kind of
examples. Partial results on 4-manifolds are also obtained.

The main result in this work is the following:

\begin{teor}
 \label{teo1} Let $f:M^n\rightarrow M^n$ be a quasi Anosov
diffeomorphism which is not Anosov. If $n=3$ then $\pi_1 (M)$
contains a subgroup which is isomorphic to ${\mathbb{Z}}^6$. If
$n=4$, $\pi_1(M)$ contains a subgroup which is isomorphic to
${\mathbb{Z}}^4$.
\end{teor}

 This result is strongly based on a theorem due to R. Plykin,
which is stated below, and
 arguments concerning Axiom A diffeomorphisms with the no-cycle condition.

\begin{teo}[\cite{ply1,ply2}]\label{teo2}
 Let $f:M\rightarrow M$ be a diffeomorphism on an
$n$-manifold , where  $n\geq 3$. If $f$ has $k$ expanding
attractors or shrinking repellers of codimension $1$ then
$\pi_1(M)$ contains a subgroup isomorphic to
$\bigoplus_k{\mathbb{Z}}^n$.
\end{teo}

We recall that a hyperbolic attractor is said to have {\em
dimension $u$} or equivalently {\em codimension $n-u$} if the
dimension of its unstable fibre bundle is $u$. Analogously we
shall say that a hyperbolic repeller is of {\em dimension $s$} or
{\em codimension $n-s$} if the dimension of its stable fiber
bundle is $s$.\par

Our strategy is to see that a QAD on a 3-manifold must have at
least two codimension one attractors or repellers, while on a
4-manifold, it must have at least one. This fact will allow us to
prove our main theorem, getting the following as an immediate
corollary:

\begin{cor}\label{conclusion.2}
Every quasi Anosov diffeomorphism on ${\mathbb T}^3$ is Anosov.
\end{cor}

Corollary \ref{conclusion.2} provides  an easier way of
recognizing Anosov diffeomorphisms on ${\mathbb T}^3$: it suffices
to check that each non zero vector in the tangent bundle goes to
infinity under the action of $Df^n$.

\section{Relation among dimension 1 attractors and codimension 1 repellers}

From now on, we shall assume $M$ is a compact smooth $n$-manifold
without boundary, and $f:M\rightarrow M$ is a QAD. The non-
wandering set of $f$ will be denoted by $\Omega(f)$, and
$\Lambda_i$ will stand for the basic sets arising from the
Spectral Decomposition Theorem \cite{sm}. We shall also maintain
the standard notation $W^s(x)$ for the stable manifold of $x$,
i.e. the set of points $y\in M$ such that $d(f^n(x),f^n(y))\to 0$
as $n\to\infty$.
\begin{defi}
Letting $\Omega(f)=\Lambda_1\cup\ldots\cup \Lambda_k$, we shall
say that $\Lambda_i$ is of type $(s,u)$ if the stable dimension of
$\Lambda_i$ is $s$, and the unstable dimension of $\Lambda_i$ is
$u$. The expression $W^\sigma(\Lambda_i)$ will be used to denote
the set $\bigcup_{x\in\Lambda_i}W^\sigma(x)$, for $\sigma=s,u$, in
which case the relation $\Lambda_i\prec\Lambda_j$ will mean that
$W^u(\Lambda_i)\cap W^s(\Lambda_j)\neq\emptyset$.
\end{defi}
We recall from \cite{m1} that a QAD is Anosov, if it satisfies the
strong transversality condition. Notice that this condition is
satisfied if all basic sets are of the same type. We state this
conclusion as a lemma

\begin{lema}
If all basic sets of a QAD $f$ are of the same type then $f$ is
Anosov.
\end{lema}

Therefore we get the following

\begin{propo}\label{prop}
Let $f$ be a QAD such that $\Omega(f)\not=M$. If some basic set of
$f$ is of type $(1,n-1)$ (or $(n-1,1)$), then $f$ must have a
codimension one attractor (repeller).
\end{propo}

\begin{proof}\,\, Let $\Lambda_1$ be a basic set of $f$, of type
$(1,n-1)$ then $W^u(\Lambda_1)$ must meet $W^s(\Lambda_2)$ where
$\Lambda_2$ is another basic set of $f$ (the no cycle condition
yields $\Lambda_1\not=\Lambda_2$). Hence $\Lambda_1\prec
\Lambda_2$, implying that the stable dimension of $\Lambda_2$ is
one. Indeed, let $x\in \Lambda_1$ and $y\in \Lambda_2$ such that
$z\in W^s(x)\cap W^u(y)$. If the stable dimension of $\Lambda_2$
were greater than one, then we would find a non zero vector $v\in
T_z M$ which would be tangent both to $W^u(z)$ and $W^s(z)$, what
would yield $\|f^n(z)v\|\to 0$ as $|n|\to \infty$, contradicting
the fact that $f$ is quasi Anosov. \par This leaves us only two
possibilities: either $W^u(\Lambda_2)$ meets $W^s(\Lambda_3)$, for
some other basic set $\Lambda_3$, or $\Lambda_2$ contains the
whole set $W^u(\Lambda_2)$, whence it would be a codimension one
attractor. In this way we inductively obtain a chain
$$\Lambda_1\prec\Lambda_2\prec\ldots\Lambda_r$$ which must end
after a finite number of steps, due to the no cycle condition.
Besides, the previous argument shows that $\Lambda_i$ has stable
dimension one, for each $i=1,\ldots r$.
\par If we suppose the previous chain is maximal, then $\Lambda_r$
must be a codimension one attractor. \end{proof}
\section{Proof of the Theorem}
We shall see that each QAD on a 3-manifold must have an attractor
{\em and} a repeller of codimension one, provided that it is not
Anosov. This, together with Theorem \ref{teo2} finishes the case
$\dim M=3$.

Let $f:M^3\rightarrow M^3$, since quasi Anosov maps do not have
attracting or repelling periodic points, their basic sets must be
all of type $(1,2)$ or $(2,1)$. If $f$ has a basic set of type
$(1,2)$ then Proposition \ref{prop} implies the existence of a
codimension one attractor. Now, if $f$ is not Anosov, then it must
also have a basic set of type $(2,1)$, otherwise, all the basic
sets would be of the same type. But in this case, Proposition
\ref{prop} guarantees the existence of a codimension one repeller,
as well.

\medskip

Observe that the previous argument also shows that if
$f:M^4\rightarrow M^4$ is not an Anosov diffeomorphism, then it
must have (at least) a basic set of type $(1,3)$ or $(3,1)$,
otherwise, all basic sets would be of type $(2,2)$ what would
imply $f$ is Anosov. Proposition \ref{prop} again, implies the
existence of an attractor or a repeller of codimension one.

\end{document}